\documentclass[12pt,letterpaper]{amsart}
\usepackage{graphicx} 
\usepackage{color}
\usepackage{amssymb}
\usepackage{MnSymbol}
\usepackage{mathtools}
\usepackage[colorlinks=true, pdfstartview=FitV, linkcolor=blue, citecolor=blue, urlcolor=blue,pagebackref=false]{hyperref}
\usepackage{microtype}
\usepackage{amsmath}
\usepackage{amsthm}
\usepackage{enumitem}
\usepackage{xifthen}
\usepackage{verbatim}
\definecolor{darkgreen}{rgb}{0,0.5,0}
\definecolor{darkblue}{rgb}{0,0,0.7}
\definecolor{darkred}{rgb}{0.9,0.1,0.1}
\usepackage[all]{xy}

\newtheorem{theorem}{Theorem}[section]

\newtheorem{prop}[theorem]{Proposition} 
\newtheorem{lemma}[theorem]{Lemma}

\theoremstyle{definition}
\newtheorem{definition}{Definition}[section]
\newtheorem{remark}[definition]{Remark}

\newtheorem{deff}[definition]{Definition}

\newcommand{\cref}[1]{Corollary~\ref{c.#1}}

\numberwithin{equation}{section}
\numberwithin{theorem}{section}

\newcommand{\C}{\mathbb{C}}

\renewcommand{\subset}{\subseteq}

\newcommand{\test}[1][]{%
\ifthenelse{\equal{#1}{}}{omitted}{given}%
}

\newcommand{\pa}{\partial}
\newcommand{\delbar}{\bar{\pa}}

\DeclareMathOperator{\Hol}{Hol}

\renewcommand{\bar}{\overline}
\renewcommand{\tilde}{\widetilde}

\renewcommand{\part}{\partial}

\usepackage{dsfont}
\usepackage{slashed}

\newcommand{\SL}{\operatorname{SL}}
\newcommand{\Tr}{\operatorname{Tr}}

\newcommand{\End}{\operatorname{End}}
\newcommand{\Hom}{\operatorname{Hom}}

\title{On a conjecture of Simpson}
\author{Panagiotis Dimakis and Sebastian Schulz}
\date{}

\begin{document}

\maketitle

\begin{abstract}
On a compact Riemann surface $\Sigma$ of genus $g\ge 2$, equipped with a complex vector bundle $E$ of rank $2$ and degree zero let $M_H$ be the moduli space of Higgs bundles. $M_H$ admits a $\mathbb C^{\star}$-action and to each stable $\mathbb C^{\star}$-fixed point $[(\bar\pa_0,\Phi_0)]$ is associated a holomorphic Lagrangian submanifold $W^1(\bar\pa_0,\Phi_0)$ inside the de Rham moduli space $M_{dR}$ of complex flat connections. In this note we prove a conjecture of Simpson \cite{Simpson} stating that $W^1(\bar\pa_0,\Phi_0)$ is closed inside $M_{dR}$. 
\end{abstract}

\tableofcontents

\section{Introduction}

Let $\Sigma$ be a compact Riemann surface of genus $g\ge 2$, equipped with a complex vector bundle $E$ of rank $n$ and degree zero. Let $M_{\mathrm{Hod}}$ be the Hodge moduli space of polystable $\lambda$-connections on $E$ up to isomorphism. $M_{\mathrm{Hod}}$ has a natural $\mathbb C^{\star}$-action and comes equipped with a canonical $\mathbb C^{\star}$-equivariant map 
\begin{equation}
    \pi: M_{\mathrm{Hod}}\to \mathbb C.
\end{equation}
The fiber $\pi^{-1}(0)$ is the moduli space $M_{\mathrm{H}}$ of polystable Higgs bundles while the fiber $\pi^{-1}(1)$ can be identified with the de Rham moduli space $M_{\mathrm{dR}}$ of complex flat connections on $E$. Every fiber is a holomorphic symplectic manifold, and indeed part of a single hyperk{\"a}hler structure. The fixed points of the $\mathbb C^{\star}$-action on $M_{\mathrm{Hod}}$ lie in $M_{\mathrm{H}}$ and will be denoted by $[(\bar\pa_0,\Phi_0)]\in M_H^{\mathbb C^{\star}}$. In this paper we only consider \textit{stable} fixed points. 

Given $p_0= [(\bar\pa_0,\Phi_0)]\in M_{\mathrm{H}}^{\mathbb C^{\star}}$, let $W(p_0)$ be the set of all points in $M_{\mathrm{Hod}}$ whose $\mathbb C^{\star}$-limit is $p_0$. The restrictions $W^0(p_0):= W(p_0)\cap\pi^{-1}(0)$ and $W^1(p_0):= W(p_0)\cap\pi^{-1}(1)$ are holomorphic Lagrangian submanifolds of $M_{\mathrm{H}}$ and $M_{\mathrm{dR}}$ respectively \cite{Simpson,CW}. 

The holomorphic Lagrangians $W^0(p_0)$ are not always closed in the ambient topology of $M_{\mathrm{H}}$. This is due to the existence of \textit{wobbly} fixed points, whose associated leaves intersect the nilpotent cone inside $M_{\mathrm{H}}$ non-trivially. See e.g.~\cite{PPN} for the case of wobbly vector bundles.

On the contrary, Simpson \cite{Simpson} conjectured that $W^1(p_0)$ are always closed inside $M_{\mathrm{dR}}$. The aim of this paper is to prove Simpson's conjecture when the rank of $E$ is $2$. 

\begin{theorem}[Main Theorem]\label{Main Theorem}
Assume that the rank of the bundle $E$ is $2$. Then for any stable $p_0= [(\bar\pa_0, \Phi_0)]\in M_{\mathrm{H}}^{\mathbb C^{\star}}$, the holomorphic Lagrangian $W^1(p_0)\subset M_{\mathrm{dR}}$ is a closed submanifold. 
\end{theorem}

The conjecture as stated is part of a larger conjecture known as the \textit{Foliation Conjecture}, which states that the individual leaves form a regular foliation of $M_{\mathrm{dR}}$. Progress has been slow, and meaningful results exist only in the case of parabolic bundles of rank two on $\mathbb{CP}^1$ minus four \cite{LSS} or five points \cite{hu2022simpson, fassarella2023moduli}. In each case, the moduli space is low-dimensional and the results are obtained via explicit computation. 

Let us briefly comment on the elements of the proof. The first element is the existence of the $\hbar$-conformal limit map $\mathcal {CL}_{\hbar}$ for $\hbar \in \C^{\star}$ \cite{DFKMMN,CW}. Given a stable $p_0 = [(\bar\pa_0, \Phi_0)]\in M_{\mathrm{H}}^{\mathbb C^{\star}}$, the map $\mathcal {CL}_{\hbar}$ descends to a biholomorphism from $W^0(p_0)$ to $W^1(p_0)$. Fixing a representative $(\bar\pa_0,\Phi_0)$ of $p_0$, the existence of a Bia\l{}ynicki-Birula slice $S^+(\bar\pa_0,\Phi_0)$ proven in \cite{CW} provides a unique distinguished representative $(\bar\pa_E,\Phi)$ for each equivalence class $p = [(\bar\pa_E,\Phi)]\in W^0(p_0)$ for which the conformal limit has the following explicit form:
\begin{equation}
    \mathcal {CL}_{\hbar}(\bar\pa_E,\Phi) = \hbar^{-1}\Phi + \bar\pa_E + \pa_0 + \hbar\Phi_0^{\dagger_{h_0}},
\end{equation}
where $h_0$ is the harmonic metric associated to $p_0$. In our context, the importance of this map is that it provides semi-explicit representatives for points in $W^1(p_0)$ in terms of points in $S^+(\bar\pa_0,\Phi_0)$. 

The second element is the proof strategy developed in \cite{Schulz}. If $W^1(p_0)$ is not closed, any point $[\nabla_{\infty}]$ in its closure can be realized as 
\begin{equation}
    [\nabla_{\infty}] = \lim_{\hbar \to 0} [\nabla_{\hbar}] =\lim_{\hbar \to 0} [\mathcal {CL}_{\hbar}(\bar\pa_E,\Phi)] = \lim_{\hbar \to 0} \left[ \hbar^{-1} \Phi + \bar\pa_E + \pa_0 + \hbar \, \Phi_0^{\dagger_{h_0}} \right],
\end{equation}
where $\hbar$ is taken to be positive real. To the family $\nabla_{\hbar}$ one can associate the trace of its \textit{holonomy} along a particular loop $\gamma$, a gauge invariant quantity. Assuming that $\Phi$ is not nilpotent, the $\mathrm{WKB}$ method developed rigorously by Mochizuki in \cite{Mochizuki, BDHH21} and the existence of a WKB loop $\gamma_{WKB}$ proven in the appendix of \cite{Schulz} show that the holonomy of $\nabla_{\hbar}$ along $\gamma_{\mathrm{WKB}}$ grows exponentially in $\hbar^{-1}$, thus proving that there can be no limit point in this case. When $\Phi$ is nilpotent, it was observed by the second author that there exists a gauge $g(\hbar)$ such that the family $\nabla'_{\hbar} = g(\hbar)\cdot \nabla_{\hbar^2}$ satisfies 
\begin{equation}
    \nabla'_{\hbar} = \hbar^{-1} \Phi' + \bar\pa_E' + \pa_E' + \mathcal O(\hbar)
\end{equation}
for some stable Higgs pair $[(\bar\pa_E',\Phi')]$ associated to $[(\bar\pa_E,\Phi)]$. 
The last element, which is the focus of this paper, is the proof that this \textit{secondary} Higgs field $\Phi'$ cannot be nilpotent unless $[(\bar\pa_E,\Phi)] = [(\bar\pa_0,\Phi_0)]$. Applying WKB to $\nabla_{\hbar}'$ thus proves Theorem~\ref{Main Theorem}. By way of contradiction we assume $\Phi'$ is nilpotent. This implies that all traces of ($n$-term products of) monodromies of $\nabla'_{\hbar}$ along non-nullhomotopic loops $\gamma:S^1\to \Sigma$ stay bounded, not only as $\hbar$ approaches zero along the positive real axis but as $|\hbar|\to 0$ along any ray. Equivalently, all traces stay bounded as $|\hbar^{-1}|\to \infty$  along any ray. Working in the original Bia\l{}ynicki-Birula gauge, we show that these traces depend holomorphically of the scaling parameter $\xi$ of the $\C^{\star}$-action. Since these traces generate the coordinate ring of the character variety and therefore separate stable points, at least one of them is a non-constant holomorphic function and therefore there exists a scaling direction along which it has to become unbounded. Since traces are gauge invariant, this contradicts the above claim and implies that $\Phi'$ cannot be nilpotent. 
\begin{remark}
The proof in this paper is completely different from the authors' previous proof attempt which contained a serious error. The proof idea in this manuscript is inspired by recent work \cite{CD} of the first author and Sotiria Chatzimarkou generalizing the structural results of \cite{CW} to Nakajima quiver varieties and using them to prove an analogue of Simpson's conjecture for arbitrary quiver varieties under a natural genericity assumption. 
\end{remark}

\subsection{Acknowledgements}
The authors would like to thank Sotiria Chatzimarkou and Georgios Kydonakis for helpful conversations and for encouragement after an error was found in a previous version of the paper. The authors thank an anonymous referee for finding the error in the previous version of this paper.

\section{Background}

\subsection{Deligne's $\lambda$-connections}

\begin{deff}
Given $\lambda \in \C$, a \textit{$\lambda$-connection} on $E$ is a triple $(\lambda, \delbar_E , \nabla_\lambda)$ where $\delbar_E$ is a Dolbeault operator on $E$ (i.e.\ the pair $(E,\delbar_E)$ is a holomorphic vector bundle $\mathcal E$) and $\nabla_\lambda : \Omega^0 (E) \rightarrow \Omega^{1,0}(E)$ is a differential operator, subject to
\begin{enumerate}[label=(\roman*)]
    \item a $\lambda$-twisted Leibniz rule: for all $f \in C^\infty (\Sigma)$ and $s \in \Omega^0 (E)$ \begin{equation}
        \nabla_\lambda (fs) = \lambda \partial f \otimes s + f \nabla_\lambda(s),
    \end{equation} 
    \item a holomorphicity condition $[ \delbar_E , \nabla_\lambda ] = 0$.
\end{enumerate}
\end{deff}

We want to study the moduli space of $\lambda$-connections up to isomorphism, meaning conjugacy-equivalence with respect to the conjugation action of the complex gauge group $\mathcal G_{\mathbb C}:= \SL(n,E)$. In order to obtain a moduli space with nice geometric properties we need to restrict to $\lambda$-connections satisfying a stability condition.

\begin{definition}
A $\lambda$- connection $(\lambda, \bar\pa_E, \nabla_{\lambda})$ is 
\begin{itemize}
    \item \textit{stable} if for any $\nabla_{\lambda}$-invariant subbundle $\mathcal F\subset \mathcal E$ it holds that $\deg \mathcal F < 0$,
    \item \textit{polystable} if it is a direct sum of stable $\lambda$-connections.
\end{itemize}
\end{definition}

We denote the space of stable (respectively polystable) $\lambda$-connections by $\mathcal M^{s}$ (respectively $\mathcal M^{ps}$).

\begin{definition}
    Define the Hodge moduli space $M_{\mathrm{Hod}}:= \mathcal M^{ps}/\mathcal G_{\mathbb C}$. 
\end{definition}

The Hodge moduli space admits a natural $\mathbb C^{\star}$-action  
\begin{equation}
    \xi\cdot [(\lambda,\bar\pa_E,\nabla_{\lambda})]:= [(\xi\lambda,\bar\pa_E,\xi\nabla_{\lambda})],
\end{equation}
as well as a $\C ^\star$-equivariant holomorphic map to $\C$ given by 
\begin{equation}
  \pi: [(\lambda, \bar\pa_E,\nabla_{\lambda})]\to \lambda .  
\end{equation}

\begin{remark}

  \begin{enumerate}[label=(\roman*)] 
    \item When $\lambda = 0$, $\Phi := \nabla_0$ becomes a holomorphic $C^\infty$-linear endomorphism of $(E, \delbar_E)$ with values in $(1,0)$-forms. Such pairs $(\delbar_E , \Phi)$ are known as \textbf{Higgs bundles} and their moduli space is denoted as $M_{\mathrm{H}}$. 
    \item When $\lambda = 1$, $D:= \bar\pa_E + \nabla_1$ is a complex flat connection and their moduli space, called the \textit{de Rham moduli space} is denoted as $M_{\mathrm{dR}}$.
    \item When $\lambda \neq 0$, then $\lambda^{-1}.[(\lambda,\bar\pa_E,\nabla_{\lambda})] \in M_{\mathrm{dR}}$. In particular, $\pi ^{-1} (\lambda) \simeq M_{\mathrm{dR}}$.
  \end{enumerate}
\end{remark}

Note that $M_{\mathrm{H}}^s:= \pi^{-1}(0)\cap\mathcal M^s/\mathcal G_{\mathbb C}$ and $M_{\mathrm{dR}}^s := \pi^{-1}(1)\cap\mathcal M^s/\mathcal G_{\mathbb C}$ are holomorphic symplectic manifolds and that they are diffeomorphic. In fact, there is a natural Hyperk{\"a}hler structure on $M_{\mathrm{H}}$ and the Hodge moduli space $M_{\mathrm{Hod}}$ is the restriction of its twistor space to the complex plane $\C \subset \mathbb{CP}^1$.

It is important to understand the fixed points of the $\mathbb C^{\star}$-action on $M_{\mathrm{Hod}}$. Clearly such fixed points must have $\lambda = 0$ and therefore lie inside the moduli space $M_{\mathrm{H}}$ of polystable Higgs bundles. 
Let $M_{\mathrm{H}}^{\mathbb C ^{\star}} \subset M_{\mathrm{H}}$ denote the set of $\mathbb C^{\star}$-fixed points $[(\bar\pa_0, \Phi_0)]$. 

\begin{prop}[Proposition $2.7$ in \cite{CW}]\label{fixed points}
    A point $[(\bar\pa_E,\Phi)]\in M_{\mathrm{H}}$ is a fixed point of the $\mathbb C^{\star}$-action if and only if there exists a splitting $E = E_1\oplus...\oplus E_l$ with respect to which 
    \begin{equation}\label{FP}
        \bar\pa_E = \begin{pmatrix} \bar\pa_{E_1} & & \\ & \ddots & \\ & & \bar\pa_{E_l} \end{pmatrix} ~~\text{and}~~ \Phi = \begin{pmatrix}0 & & & \\ \Phi_1 & 0 & & \\ & \ddots & \ddots & \\ & & \Phi_{l-1} & 0 \end{pmatrix}
    \end{equation}
    where $\Phi_j : E_j \to E_{j+1}\otimes K$ is holomorphic. 
\end{prop}

Let us also define the following endomorphism spaces which will be used in the next subsection:
\begin{equation}
    \begin{split}
        N_+ :=& \End_{>0} E = \bigoplus\limits_{i-k>0} \Hom (E_i,E_k),\\
        L:=& \left(\bigoplus\End(E_i)\right)\cap \mathfrak{sl}(E).
    \end{split}
\end{equation}

\subsection{Holomorphic Lagrangians and the conformal limit}

Given a \textit{stable} Higgs pair $[(\bar\pa_E,\Phi)]$ let $p_0= [(\bar\pa_0,\Phi_0)]:=\lim\limits_{\xi\to 0}\xi\cdot [(\bar\pa_E,\Phi)]$ be its $\mathbb C^{\star}$-limit. Define 
\begin{equation}
    W(p_0):= \{ [(\lambda, \bar\pa_E,\nabla_{\lambda})]\in M_{\mathrm{Hod}}: \lim\limits_{\xi\to 0} \xi\cdot [(\lambda, \bar\pa_E,\nabla_{\lambda})] = p_0\}.
\end{equation}
Let $W^i(p_0):= W(p_0)\cap\pi^{-1}(i)$ for $i\in\{0,1\}$. It is proven in \cite{Simpson,CW} that the submanifolds $W^i(p_0)$ are the leaves of holomorphic Lagrangian foliations of $M_{\mathrm{H}}$ and $M_{\mathrm{dR}}$ respectively. These holomorphic Lagrangian sub-manifolds interact nicely with the conformal limit construction of \cite{CW}. In order to explain the relation let us fix a representative  $(\bar\pa_0,\Phi_0)$ of a $\C^{\star}$-fixed point $p_0$.  
\begin{definition}[\cite{CW}]
The Bia\l{}ynicki-Birula slice $S^+(\bar\pa_0,\Phi_0)$ is defined as 
\begin{equation}
\begin{split}
S^+(\bar\pa_0,\Phi_0):= \{(\beta,\phi)\in \Omega^{0,1}(N_+)\oplus\Omega^{1,0}(L\oplus N_+), \\
D''(\beta,\phi) + [\beta,\phi]= 0,~~D'(\beta,\phi) = 0\},
\end{split}
\end{equation}
where 
\begin{equation}
    D'':= \bar\pa_0 + \Phi_0,~~D':= \pa_0^{\dagger_{h_0}} + \Phi_0^{\dagger_{h_0}}.
\end{equation}
\end{definition}
Here, $h_0$ is the \textit{harmonic metric} of the Higgs bundle $(\delbar_0, \Phi_0)$. That is, it is the unique hermitian metric on $E$ which solves the \textit{Hitchin equation}
\begin{equation}
    F_{D_0} + \left[ \Phi_0, \Phi_0^{\dagger_{h_0}} \right] = 0,
\end{equation}
where $D_0 = \delbar_0 + \pa _0^{\dagger_{h_0}}$ is the Chern connection of the pair $(\delbar_0, h_0)$.

It is proven in \cite{CW} that $S^+(\bar\pa_0,\Phi_0)$ not only is biholomorphic to $W^0(p_0)$ but that each element $[(\bar\pa_E,\Phi)]\in W^0(p_0)$ is uniquely gauge equivalent to $(\bar\pa_0+\beta,\Phi_0+\phi)$ with $(\beta,\phi)\in S^+(\bar\pa_0,\Phi_0)$. We say that a representative of $W^0(p_0)$ is in BB-gauge if it is in the above form. It is then proven that there exists a biholomorphic map $\mathcal {CL}_{\hbar}: W^0(p_0)\to W^1(p_0)$ called the conformal limit, which takes the following explicit form when the point is put in BB-gauge:
\begin{equation}
    \mathcal {CL}_{\hbar}(\bar\pa_0+\beta,\Phi_0+\phi) = \hbar^{-1}(\Phi_0+\phi) + \bar\pa_0 +\beta + \pa_0^{\dagger_{h_0}} + \hbar\Phi_0^{\dagger_{h_0}}.
\end{equation}
\begin{lemma}
If the representative of $[(\bar\pa_E,\Phi)]$ in BB gauge is given by 
\begin{equation*}
(\bar\pa_0 + \sum\limits_{j=1}^{l-1} \beta_j, \Phi_0 + \sum\limits_{j=0}^{l-1} \phi_j) 
\end{equation*}
with $\beta_j \in \Omega^{0,1}\left(\bigoplus\limits_{i-k=j}\mathrm{Hom}(E_i,E_k)\right)$ and $\phi_j \in \Omega^{1,0}\left(\bigoplus\limits_{i-k=j}\mathrm{Hom}(E_i,E_k)\right)$ then the representative of $\xi\cdot[(\bar\pa_E,\Phi)]$ is given by 
\begin{equation*}
\left(\bar\pa_0 + \sum\limits_{j=1}^{l-1} \xi^j\beta_j, \Phi_0 + \sum\limits_{j=0}^{l-1} \xi^{j+1}\phi_j\right) =: (\bar\pa_0+\beta(\xi),\Phi_0+\phi(\xi)).
\end{equation*}
\end{lemma}
The proof of this lemma is straightforward. Since the harmonic metric $h_0$ respects the grading at the fixed point, the above lemma implies that for the gauge transformation 
\begin{equation}
g(\xi) := \mathrm{diag}(\xi^{-\frac{l-1}{2}}\mathbf{I}_{E_1}, \xi^{-\frac{l-3}{2}}\mathbf{I}_{E_2},...,\xi^{\frac{l-1}{2}}\mathbf{I}_{E_l}),
\end{equation}
where $\mathbf{I}_{E_j}$ denotes the identity endomorphism restricted to $E_j$, 
\begin{equation}
\begin{split}
g(\xi)\cdot \mathcal {CL}_{\hbar}(\bar\pa_0+\beta(\xi),\Phi_0+\phi(\xi)) = \hbar^{-1}\xi(\Phi_0+\phi) + \bar\pa_0 +\beta + \pa_0^{\dagger_{h_0}} + \hbar\xi^{-1}\Phi_0^{\dagger_{h_0}}.
\end{split}
\end{equation}

\section{Proof of Simpson's conjecture}

\subsection{WKB and nilpotent Higgs fields}

\cite{Schulz} contains partial results and a proof strategy for Simpson's conjecture which we review now.

Assume that for some fixed point $p_0$, the boundary $B:=\pa W^1(p_0)$ of the leaf $W^1(\bar\pa_0, \Phi_0)\subset M_{dR}$ is not closed and fix $\nabla\in B$.

The first observation is that $\nabla$ can be realized as the limit of a family of complex flat connections coming from the Conformal Limit construction. To be precise, there exists $(\beta,\phi) \in S^+(\bar\pa_0, \Phi_0)$ such that
\begin{equation}\label{boundary point}
\begin{split}
    [\nabla] &= \lim_{\xi \to +\infty} [g(\xi)\cdot \mathcal {CL}_{\hbar}(\bar\pa_0+\beta(\xi),\Phi_0+\phi(\xi))] \\
    &=\lim_{\xi \to +\infty} [ \hbar^{-1} \xi \Phi + \bar\pa_E + \pa_0^{\dagger_{h_0}} + \hbar\xi^{-1} \, \Phi_0^{\dagger_{h_0}}],
\end{split}
\end{equation}
where $(\bar\pa_E,\Phi) = (\bar\pa_0+\beta,\Phi_0+\phi)$. Here, the limit is taken along the positive real axis.

To arrive at a contradiction we want to show that the right hand side cannot have a finite limit point. To achieve this we want to associate to it a gauge-invariant quantity which diverges as $\xi\to \infty$. For ease of notation, we rename $\hbar^{-1}\xi =r$ and we want to investigate when does the curve comprised of equivalence classes of the family 
\begin{equation}
\nabla_r:= r \Phi + \bar\pa_E + \pa_0^{\dagger_{h_0}} + r^{-1} \, \Phi_0^{\dagger_{h_0}}
\end{equation}
in $M_{\mathrm{dR}}$ becomes unbounded as $r\to \infty$. The quantity we use is the trace of the holonomy around some particular closed loop on the surface. The method used to study the asymptotic behavior of such quantities is known as the \textit{WKB method}. In order to state the main result of the WKB method we need to introduce some terminology. 

Fix a rank $2$ holomorphic vector bundle $\mathcal E$ and a family of complex flat connections 
\begin{equation}\label{r family}
    \nabla_{r} = r^{-1}\varphi + D + r\psi.
\end{equation}
Let $\gamma: S^1_t\to \Sigma$ be a parametrized loop and denote by $A$ the $1$-form that represents $\gamma^{\star}D$. Assuming the Higgs field $\varphi$ is \textbf{not} nilpotent, the pullback $\gamma^{\star}\varphi$ has two eigenvalues $\pm \mu(t)\,dt.$ We say that $\gamma$ is a \textit{WKB curve} for $\varphi$ if $\Re(\mu(t))>0$ $\forall t$. For any $\varphi$ there exist closed loops on $\Sigma$ which are WKB curves for $\varphi$. A proof of this fact can be found in the appendix of \cite{Schulz}. From this point on we assume that $\gamma$ is indeed a WKB curve. Let $Z_{\gamma}:= \int\limits_{S^1_t} \mu(t)\,dt$. Clearly, along a WKB curve $\Re (Z_{\gamma})>0$. Since $\pm\mu(t)$ are distinct, we can choose a splitting $\gamma^{\star}\mathcal E= L_{-}\oplus L_{+}$ which diagonalizes $\gamma^{\star}\varphi$. With respect to this decomposition we write $A= A_{+}+A_{-} + A_0$ with $A_{\pm}$
denoting the $1$-forms representing flat metrics on $L_{\pm}$ and $A_0$ being the off-diagonal part of the connection. The WKB method then provides the following: 

\begin{theorem}\label{WKB}
    Along a WKB curve $\gamma$,
    \begin{equation}
        \lim\limits_{r\to \infty}(\Tr\Hol_{\gamma}(\nabla_r)\exp(-rZ_{\gamma}) = \Hol_{\gamma}(A_+)\in \mathbb C^{\star}.
    \end{equation}
\end{theorem}

A proof by Mochizuki can be found in the appendix of \cite{BDHH21}. Unfortunately, in the family of complex flat connections appearing in \eqref{boundary point} the Higgs field $\Phi$ is nilpotent so Theorem~\ref{WKB} does not directly apply. 

Since $\Phi$ is nilpotent, $L_1 := \ker\Phi$ is a holomorphic line sub-bundle of $\mathcal E$. Let $L_2$ be the orthogonal complement of $L_1$ in $\mathcal E$ with respect to the harmonic metric associated to $(\bar\pa_E,\Phi)$. With respect to the splitting $E= L_1\oplus L_2$ we can write 
\begin{equation}
    \bar\pa_E + \pa_0^{\dagger_{h_0}} = D = A_1 + D' + A_{-1}
\end{equation}
with $A_1\in \Gamma(\Hom(L_2,L_1)\otimes\Omega^1)$, $A_{-1}\in \Gamma(\Hom(L_1,L_2)\otimes\Omega^1)$ and $D'$ diagonal. In these coordinates 
\begin{equation}
    \Phi = \begin{pmatrix}0 & \tilde\Phi \\ 0 & 0 \end{pmatrix}.
\end{equation}
Let $g(r)$ be the constant gauge transformation
\begin{equation}
    g(r)= \begin{pmatrix}r^{-1/2} &0 \\ 0 & r^{1/2} \end{pmatrix}
\end{equation}
and consider 
\begin{equation}
    \nabla_r':= g(r)\cdot \nabla_{r^2} = r(\Phi + A_{-1}) + D' + \mathcal O(r^{-1}).
\end{equation}
If $\bar\pa_E'$ is the $(0,1)$ part of $D'$ and $\Phi' = \Phi+ A_{-1}$, then Proposition $3.2$ of \cite{Schulz} shows that the pair $(\bar\pa_E',\Phi')$ is a stable Higgs bundle. Therefore, in order to be able to apply Theorem~\ref{WKB} to the family of connections in \eqref{boundary point} we need to show that the Higgs field $\Phi'$ is not nilpotent. Notice that with respect to the splitting $E= L_1\oplus L_2$, 
\begin{equation}
    \Phi' = \begin{pmatrix} 0 & \tilde\Phi \\ A_{-1} & 0 \end{pmatrix},
\end{equation}
so we need to show that $A_{-1}\neq 0$. 

\subsection{Proof of the main theorem}

By way of contradiction we assume that $A_{-1} =0$. This implies that the family $\nabla_r$ can be written in the form 
\begin{equation}
\nabla_r = d + r\begin{pmatrix}0 & \tilde\Phi \\ 0 & 0 \end{pmatrix} + \begin{pmatrix} d_1 & A_1 \\ 0 & d_2 \end{pmatrix} + r^{-1}M
\end{equation}
with respect to the splitting $E = L_1\oplus L_2$. All matrices are valued in $\Omega^{1}$ and all entries $d_1,d_2,A_1$ and the entries of $M$ are bounded. Moreover, this form is valid not only for positive real $r$ but in fact for all $r\in \C^{\star}$. Let us write explicitly 
\begin{equation}
M = \begin{pmatrix} m_{11} & m_{12} \\ m_{21} & m_{22} \end{pmatrix} 
\end{equation}
and consider 
\begin{equation}
g(r)\cdot \nabla_r = d + \begin{pmatrix} d_1 + r^{-1}m_{11} & \tilde\Phi + m_{12} \\ m_{21} & d_2 + r^{-1}m_{22}\end{pmatrix}.
\end{equation}
When $|r|$ stays uniformly away from zero all coefficients are uniformly bounded and therefore given any finite set of non-nullhomotopic loops $\gamma_i$ in $\Sigma$ (possibly avoiding the finitely many points where $L_1$ and $L_2$ may become dependent,) the traces 
\begin{equation}
\Tr(\Hol_{\gamma_1}(\nabla_r)\cdot\Hol_{\gamma_2}(\nabla_r)\cdot...\cdot\Hol_{\gamma_n}(\nabla_r)) 
\end{equation}
stay bounded as $|r|\to \infty$. We now proceed to show that this cannot happen for all $|r|\to \infty$. Before we give the final argument we need the following theorem for the $\SL(2,\C)$-character variety $X_{\SL(2,\C)}(\pi_1(\Sigma))$.
\begin{theorem}[Proposition $2$ in \cite{Sikora}]\label{Sikora}
Let $\gamma_1,...,\gamma_{2g}$ be loops representing the generators of $\pi_1(\Sigma)$. The algebra $\C[X_{\SL(2,\C)}(\pi_1(\Sigma))]$ is generated by the elements $\Tr(\Hol_{\gamma_i}(\nabla))$, $\Tr(\Hol_{\gamma_i}(\nabla)\cdot\Hol_{\gamma_j}(\nabla))$ and $\Tr(\Hol_{\gamma_i}(\nabla)\cdot\Hol_{\gamma_j}(\nabla)\cdot\Hol_{\gamma_k}(\nabla))$ for $i<j<k$. 
\end{theorem}
It follows from the GIT construction of the character variety that these functions separate closed orbits. Going back to the BB gauge representation of our family of connections, we can write them as 
\begin{equation}
\nabla_{\xi} = \bar\pa_0 + \sum\limits_{j=1}^{l-1} \xi^j\beta_j +\hbar^{-1}\left( \Phi_0 + \sum\limits_{j=0}^{l-1} \xi^{j+1}\phi_j\right) + \pa_0^{\dagger_{h_0}} + \Phi_0^{\dagger_{h_0}}. 
\end{equation}
Notice that in this description we keep $\hbar$ fixed. If we pull back this connection along a  non-nullhomotopic loop $\gamma$, it obtains the following form 
\begin{equation}
\gamma^{\star}\nabla_{\xi} = \frac{\,d}{\,dt}\,dt + (A_0(t)+ \xi A_1(t)+...+ \xi^{l} A_l(t))\,dt. 
\end{equation}
In fact for the group $\SL(2,\C)$, $l\le 2$. This in particular implies that the monodromy along the loop $\gamma$ satisfies a first order ordinary differential equation whose zero order coefficient is an analytic function on the parameter $\xi$. Standard ordinary differential equation theory implies that the solution depends holomorphically on $\xi$. In particular all the trace functions appearing in Theorem \ref{Sikora} become holomorphic functions in $\xi$. If at least one of them is not constant then standard complex analysis implies that there is a ray emanating from the origin such that this function becomes unbounded as $|\xi|\to \infty$ along this ray. In fact by the Phragm\'en--Lindel\"of principle at least half the directions work. This would contradict the conclusion we reached above for $r = \hbar^{-1}\xi$, therefore, all traces have to be constant. Note that these trace functions distinguish irreducible representations and, indeed, the representations we consider here are irreducible. This is because we know that the conformal limit maps the $\C^{\star}$-fixed point $[(\bar\pa_0,\Phi_0)]$ to its image under the Non-Abelian Hodge correspondence and therefore the associated representation is irreducible. Since the conformal limit is a biholomorphism onto its image when restricted to $W^0(p_0)$ and irreducibility is an open condition, for all $\xi$ sufficiently small the connections $\nabla_{\xi}$ have to be irreducible and distinct. Therefore at least one of the traces has to be non-constant when restricted to this family. Finally we have reached a contradiction. 

\bibliography{bibliography}{}
\bibliographystyle{amsalpha}

\smallbreak
 
\noindent {\bf  Panagiotis Dimakis - 
{\sc  University of Maryland, USA.}\\
\tt  pdimakis12345@gmail.com } \\

\noindent {\bf  Sebastian Schulz  \\
\tt  sschulzmath@gmail.com }

\end{document}